\documentclass[MSNbibl,nameyear,dvips]{arxstspdf}
\usepackage{flushend}
\usepackage{stfloats}

\volume{29}
\issue{2}
\pubyear{2014}
\firstpage{247}
\lastpage{251}
\doi{10.1214/14-STS472} 
\referstodoi{10.1214/13-STS457}

\makeatletter

\newproclaim{ep}{Efficiency principle}
\newproclaim{vp}{Validity principle}
\newproclaim{example}{Example}
\newproclaim{astep}{A-step}
\newproclaim{pstep}{P-step}
\newproclaim{cstep}{C-step}

\newcommand{\prob}{\mathsf{P}}

\newcommand{\bel}{\mathsf{bel}}
\newcommand{\pl}{\mathsf{pl}}

\newcommand{\unif}{\mathsf{Unif}}
\newcommand{\nm}{\mathsf{N}}

\newcommand{\XX}{\mathbb{X}}
\newcommand{\UU}{\mathbb{U}}

\renewcommand{\S}{\mathcal{S}}

\renewcommand{\citep}[1]{(\citeauthor{#1}, \citeyear{#1})}
\newcommand{\eqref}[1]{(\ref{#1})}
\newcommand{\st}{\mathrm{st}}

\makeatother

\begin{document}
\begin{frontmatter}

\vspace*{6pt}\title{Discussion: Foundations of Statistical Inference, Revisited}%
\runtitle{Discussion}

\begin{aug}
\author[a]{\fnms{Ryan} \snm{Martin}\corref{}\ead[label=e1]{rgmartin@uic.edu}}
\and
\author[b]{\fnms{Chuanhai} \snm{Liu}\ead[label=e2]{chuanhai@purdue.edu}}
\runauthor{R. Martin and C. Liu}

\affiliation{University of Illinois at Chicago and Purdue University}

\address[a]{Ryan Martin is Assistant Professor,
Department of Mathematics, Statistics, and Computer Science, University of Illinois at Chicago,
851 S. Morgan St., Chicago, Illinois 60607,
USA
\printead{e1}.}

\address[b]{Chuanhai Liu is Professor,
Department of Statistics, Purdue University,
250 North University St.,
West Lafayette, Indiana 47907-2067,
USA
\printead{e2}.}
\end{aug}

\begin{abstract}
This is an invited contribution to the discussion on Professor
Deborah Mayo's paper, ``On the Birnbaum argument for the strong
likelihood principle,'' to appear in \emph{Statistical Science}.
Mayo clearly demonstrates that statistical methods violating the
likelihood principle need not violate either the sufficiency or
conditionality principle, thus refuting Birnbaum's claim. With
the constraints of Birnbaum's theorem lifted, we revisit the
foundations of statistical inference, focusing on some new
foundational principles, the inferential model framework, and
connections with sufficiency and conditioning.
\end{abstract}

%
\end{frontmatter}

\section{Introduction}
\label{S:intro}

Birnbaum's theorem \citep{birnbaum1962} is arguably the most
controversial result in statistics. The theorem's conclusion is that a
framework for statistical inference that satisfies two natural
conditions, namely, the sufficiency principle (SP) and the
conditionality principle (CP), must also satisfy an exclusive
condition, the likelihood principle (LP). The controversy lies in the
fact that LP excludes all those standard methods taught in Stat 101
courses. Professor Mayo successfully refutes Birnbaum's claim, showing
that violations of LP need not imply violations of SP or CP. The key to
Mayo's argument is a correct formulation of CP; see also \citet
{evans2013}. Her demonstration resolves the controversy around Birnbaum
and LP, helping to put the statisticians' house in order.

The controversy and confusion surrounding Birnbaum's claim has perhaps
discouraged researchers from considering questions about the
foundations of statistics. We view Professor Mayo's paper as
an invitation for statisticians to revisit these fundamental questions,
and we are grateful for the opportunity to contribute to this discussion.

Though LP no longer constrains the frequentist approach, this does not
mean that pure frequentism is necessarily correct. For example,
reproducibility issues\footnote{``Announcement: Reducing our
irreproducibility,'' \emph{Nature} \textbf{496} (2013),
DOI:\doiurl{10.1038/496398a}.} in large-scale studies is an indication that the
frequentist techniques that have been successful in classical problems
may not be appropriate for today's high-dimensional problems. We
contend that something more than the basic sampling model is required
for valid statistical inference, and appropriate conditioning is one
aspect of this. Here we consider what a new framework, called \emph
{inferential models} (IMs), has to say concerning the foundations of
statistical inference, with a focus on something more fundamental than
CP and SP. For this, we begin in Section~\ref{S:valid} with a
discussion of valid probabilistic inference as motivation for the IM
framework. A general efficiency principle is presented in Section~\ref{S:drp} and we give an IM-based dimension reduction strategy that
accomplishes what the classical SP and CP set out to do. Section~\ref{S:discuss} gives some concluding remarks.

\section{Valid Probabilistic Inference}
\label{S:valid}

\subsection{A Validity Principle}

The sampling model for observable data $X$, depending on an unknown
parameter $\theta$, is the familiar starting point. Our claim is that
the sampling model alone is not sufficient for valid probabilistic
inference. By ``probabilistic inference'' we mean a framework whereby
any assertion/hypothesis about the unknown parameter has a (predictive)
probability attached to it after data is observed. The sampling model
facilitates a comparison of the chances of the event $X=x$ on different
probability spaces. For probabilistic inference, there must be a known
distribution available, after data is observed. The random element that
corresponds to this distribution shall be called a predictable
quantity, and probabilistic inference is obtained by predicting this
predictable quantity after seeing data. The probabilistic inference is
``valid'' if the predictive probabilities are suitably calibrated or,
equivalently, have a fixed and known scale for meaningful
interpretation. Without risk of confusion, we shall call the prediction
of the predictable quantity valid if it admits valid probabilistic
inference. We summarize this in the following \emph{validity
principle} (VP).
%
\begin{vp*}
Probabilistic inference requires associating an unobservable but
predictable quantity with the observable data and unknown parameter.
Probabilities to be used for inference are obtained by valid prediction
of the predictable quantity.
\end{vp*}

The frequentist approach aims at developing  procedures, such
as confidence intervals and testing rules, having long-run frequency
properties. Expressions like ``95\% confidence'' have no predictive
probability interpretation after data is observed, so frequentist
methods are not probabilistic in our sense.
Nevertheless, certain frequentist quantities, such as $p$-values, may
be justifiable from a valid probabilistic inference point of view
\citep{impval}.

When genuine prior information is available and can be summarized as a
usual probability model, the corresponding Bayesian inference is both
probabilistic and valid [see \citet{imcond}, Remark~4]. When no genuine
prior information is available, and a default prior distribution is
used, the validity property is questionable. Probability matching
priors, Bernstein--von Mises theorems, etc., are efforts to make
posterior inference valid, in the sense above, at least approximately.
The standard interpretation of these results is, for example, that
Bayesian credible intervals have the nominal frequentist coverage
probability asymptotically; see \citet{fraser2011}. In that case, the
remarks above concerning frequentist methods apply.

Fiducial inference was introduced by \citet{fisher1930} to avoid using
artificial priors in scientific inference. Subsequent work includes
structural inference \citep{fraser1968}, the Dempster--Shafer theory
of belief functions (\cite{shafer1976}; \cite{dempster2008}), generalized
inference (\cite{chiang2001}; \cite{weerahandi1993}) and generalized fiducial
inference (\citeauthor{hannig2009}, \citeyear{hannig2009}, \citeyear{hannig2012}). Fiducial distributions are
defined by expressing the parameter as a data-dependent function of a
pivotal quantity. This results in a bona fide posterior distribution
only in Fraser's structural models and, in those cases, it corresponds
to a Bayesian posterior (\cite{lindley1958}; \cite{taraldsen.lindqvist.2013}).
Therefore, the fiducial distribution is meaningful when the
corresponding Bayesian prior is meaningful in the sense above. More on
fiducial from the IM perspective is given below.

\subsection{IM Framework}
\label{SS:im.review}

The IM framework, proposed recently by \citet{imbasics}, has its roots
in fiducial and Dempster--Shafer theory; see also \citet{mzl2010}. At
a fundamental level, the IM approach is driven by VP. Here is a quick overview.

Write the sampling model/data-generating mechanism as
%
\begin{equation}
\label{eq:assoc} X = a(\theta, U),\quad  U \sim\prob_U,
\end{equation}
where $X \in\XX$ is the observable data, $\theta\in\Theta$ is the
unknown parameter, and $U \in\UU$ is an unobservable auxiliary
variable with known distribution $\prob_U$. Following VP, the goal is
to use the data $X$ and the distribution for $U$ for meaningful
probabilistic inference on $\theta$ without assuming a prior. The
following three steps describe the IM construction.
%
\begin{astep*}
Associate the observed data $X=x$, the parameter and the auxiliary
variable via \eqref{eq:assoc} and construct the set-valued mapping,
given by
\[
\Theta_x(u) = \bigl\{\theta\dvtx  x=a(\theta,u)\bigr\}, \quad u \in\UU.
\]
The fiducial approach considers the distribution of $\Theta_x(U)$ as a
function of $U \sim\prob_U$. The IM framework, on the other hand,
predicts the unobserved $U$ using a random set.
\end{astep*}
%
\begin{pstep*}
Predict the unobservable $U$ with a random set $\S$. The distribution
$\prob_\S$ of $\S$ is required to be valid in the sense that
%
\begin{equation}
\label{eq:valid} f(U) \geq_{\st } \unif(0,1),
\end{equation}
where $f(u) = \prob_\S(\S\ni u)$ and $\geq_{\st }$ means
``stochastically no smaller than.''
\end{pstep*}
%
\begin{cstep*}
Combine $\Theta_x(\cdot)$ and $\S$ as $\Theta_x(\S) = \linebreak[4] \bigcup_{u
\in\S} \Theta_x(u)$. For a given assertion $A \subseteq\Theta$,
evaluate the evidence in $x$ for and not against the claim ``$\theta
\in A$'' via the belief and plausibility functions:
\begin{eqnarray*}
\bel_x(A) & = &\prob_\S\bigl\{\Theta_x(\S)
\subseteq A\bigr\},
\\
\pl_x(A) & =& \prob_\S\bigl\{\Theta_x(\S)
\cap A \neq\varnothing\bigr\}.
\end{eqnarray*}
In cases where $\Theta_x(\S)$ is empty with positive $\prob_\S
$-probability, some adjustments to the above formulas are needed; see
\citet{imbasics}.
\end{cstep*}

The IM output is the pair of functions $(\bel_x,\pl_x)$ and, when
applied to an assertion $A$ about the parameter $\theta$ of interest,
these provide a (personal) probabilistic summary of the evidence in
data $X=x$ supporting the truthfulness of $A$. Property \eqref
{eq:valid} guarantees that the IM output is valid. Our focus in the
remainder of the discussion will be on the A-step and, in particular,
auxiliary variable dimension reduction. Such concerns about
dimensionality are essential for efficient inference on $\theta$.
These are also closely tied to the classical ideas of sufficiency and
conditionality.

We end this section with a few remarks on fiducial. If one takes the
predictive random set $\S$ in the P-step as a singleton, that is, $\S
= \{U\}$, where $U \sim\prob_U$, then $\bel_x$ and $\pl_x$ are
equal and equal to the fiducial distribution. In this sense, fiducial
provides probabilistic inference. However, the singleton predictive
random set is not valid in the sense of \eqref{eq:valid}, so fiducial
inference is generally not valid, violating the second part of VP. One
can also reconstruct the fiducial distribution by choosing a valid
predictive random set $\S$ so that $\bel_x(A)$ equals the fiducial
probability of $A$ for all suitable $A \subseteq\Theta$. But this
would generally require that $\S$ depend on the observed $X=x$, and
the resulting inference suffers from a selection bias, or double-use of
the data, resulting in unjustifiably large belief probabilities.

\section{Efficiency and Dimension Reduction}
\label{S:drp}

\subsection{An Efficiency Principle}

It is natural to strive for efficient statistical inference. In the
context of IMs, we want $\pl_X(A)$ to be as stochastically small as
possible, as a function of $X$, when the assertion $A$ about $\theta$
is false. To connect this to classical efficiency, $\pl_X(A)$ can be
interpreted like the $p$-value for testing $H_0\dvtx  \theta\in A$, so
stochastically small plausibility when $A$ is false corresponds to the
high power of the test. We state the following \emph{efficiency
principle} (EP).
%
\begin{ep*}
Subject to the validity constraint, probabilistic inference should be
made as efficient as possible.
\end{ep*}

EP is purposefully vague: it allows for a variety of techniques to be
employed to increase efficiency. The next section discusses one
important technique related to auxiliary variable dimension reduction.

\subsection{Improved Efficiency via Dimension Reduction}

In the classical setting, sufficiency reduces the data to a good
summary statistic. In the IM context, however, the dimension of the
auxiliary variable, not the data, is of primary concern. For example,
in the case of iid sampling, the dimension of $U$ is the same as that
of $X$, which is usually greater than that of $\theta$. In such cases,
it is inefficient to predict a high-dimensional auxiliary variable for
inference on a lower dimensional parameter. The idea is to reduce the
dimension of $U$ to that of $\theta$. This auxiliary variable
dimension reduction will indirectly result in some transformation of
the data.

How to reduce the dimension of $U$? We seek a new auxiliary
variable $V$, of the same dimension of $\theta$, such that the
baseline association \eqref{eq:assoc} can be rewritten as
%
\begin{equation}
\label{eq:assoc1} T(X) = b(\theta, V), \quad V \sim\prob_V,
\end{equation}
for some functions $T$ and $b$, and distribution $\prob_V$. Here
$\prob_V$ may actually depend on some features of the data $X$. Such a
dimension reduction is general, but \citet{imcond} consider an
important case, which we summarize here. Suppose we have two one-to-one
mappings, $x \mapsto(T(x),H(x))$ and $u \mapsto(\tau(u), \eta(u))$,
with the requirement that $\eta(U)=H(X)$. Since $H(X)$ is observable,
so too must be the feature $\eta(U)$ of $U$. This point has two
important consequences: first, a feature of $U$ that is observed need
not be predicted, hence a dimension reduction; second, the feature of
$U$ that is observed naturally provides some information about the part
that remains unobserved, so conditioning should help improve
prediction. By construction, the baseline association
\eqref{eq:assoc} is equivalent to
%
\begin{equation}
\label{eq:assoc2} T(X) = b\bigl(\theta, \tau(U)\bigr) \quad \mbox{and}\quad  H(X) = \eta(U),
\end{equation}
and this suggests an association of the form \eqref{eq:assoc1}, where
$V=\tau(U)$ and $\prob_V$ is the conditional distribution of $\tau
(U)$ given $\eta(U)=H(X)$.

It remains to discuss how the dimension reduction strategy described
above related to EP. The following theorem gives one relatively simple
illustration of the improved efficiency via dimension reduction.
%
\begin{thm}
\label{thm:drp}
Suppose that the baseline association \eqref{eq:assoc} can be
rewritten as \eqref{eq:assoc2} and that $\tau(U)$ and $\eta(U)$ are
independent. Then inference based on $T(X) = b(\theta,\tau(U))$
alone, by a valid prediction of $\tau(U)$, ignoring $\eta(U)$, is at
least as efficient as inference from \eqref{eq:assoc2} by a valid
prediction of $(\tau(U),\eta(U))$.
\end{thm}

See the corollary to Proposition~1 in Martin and Liu (\citeyear{imcond}, full-length
version); see also \citet{imbook}. Therefore, reducing the
dimension of the auxiliary variable cannot make inference less
efficient. The point in that paper is that reducing the dimension
actually improves efficiency, hence EP.

In the standard examples, for example, regular exponential families,
our dimension reduction above corresponds to classical sufficiency; see
Example~\ref{ex:normal}. Outside the standard examples, the IM
dimension reduction gives something different from sufficiency, in
particular, the former often leads directly to further dimension
reduction compared to the latter; see Example~\ref{ex:nile}. That the
IM dimension reduction naturally contains some form of conditioning is
an advantage. The absence of conditioning in the standard definition of
sufficiency is one possible reason why conditional inference has yet to
become part of the mainstream. The IM framework also has advantages
beyond dimension reduction and conditioning. In particular, the IM
output gives valid prior-free probabilistic inference on $\theta$.

\subsection{``Dimension Reduction Entails SP and CP''}
\label{SS:entails}

This section draws some connections between the IM dimension reduction
above and SP and CP. First, it is clear that following the auxiliary
variable dimension reduction strategy described above entails CP. In
the Cox example, the randomization that determines which measurement
instrument will be used corresponds to an auxiliary variable whose
value is observed completely. So, our auxiliary variable dimension
reduction strategy implies conditioning on the actual instrument used,
hence CP. For SP, Theorem~\ref{thm:drp} gives some insight. That is,
when a sufficient statistic has dimension the same as $\theta$, one
can take $T(X)$ as that sufficient statistic and select independent
$\tau(U)$ and $\eta(U)$. In general, our dimension reduction and
efficiency considerations are more meaningful than sufficiency and SP.
The examples below illustrate this point further.
%
\begin{example}
\label{ex:normal}
Suppose $X_1,X_2$ are independent $\nm(\theta,1)$ samples, and write
the association as $X_i=\theta+ U_i$, where $U_1,U_2$ are independent
standard normal. In this case, there are lots of candidate mappings
$(\tau, \eta)$ to rewrite the baseline association in the form \eqref
{eq:assoc2}. Two choices are $\{\tau(u)=u_1, \eta(u)=u_2-u_1\}$ and
$\{\tau(u)=\bar u, \eta(u)=(u_1-\bar u, u_2-\bar u)\}$. At first
look, the second choice, corresponding to sufficiency, seems better
than the first. However, the dimension-reduced associations \eqref
{eq:assoc1} based on these two choices are exactly the same. This
means, first, there is nothing special about sufficiency in light of
proper conditioning (\cite{evans.fraser.monette.1986}; \cite{fraser2004}).
Second, it suggests that, at least in simple problems, the
dimension-reduced association \eqref{eq:assoc1} does not depend on the
choice of $(\tau, \eta)$, that is, it only depends on the sufficient
statistic, hence SP. The message here holds more generally, though a
rigorous formulation remains to be worked out.
\end{example}
%
\begin{example}
\label{ex:nile}
Consider independent exponential random variables $X_1,X_2$, the first
with mean $\theta$ and the second with mean $\theta^{-1}$. In this
case, the minimal sufficient statistic, $(X_1,X_2)$, is two-dimensional
while the parameter is one-dimensional. \citet{imcond} take the
baseline association as $X_1=\theta U_1$ and $X_2 = \theta^{-1} U_2$,
where $U_1,U_2$ are independent standard exponential. They employ a
novel partial differential equations technique to identify a function
$\eta$ of $(U_1,U_2)$ whose value is fully observed, so that only a
scalar auxiliary variable needs to be predicted. Their solution is
equivalent to that based on the conditional distribution of the maximum
likelihood estimate given an ancillary statistic (\cite{fisher1973};
\cite{ghoshreidfraser2010}). The message here is that the IM-based auxiliary
variable dimension reduction strategy does something similar to the
classical strategy of conditioning on ancillary statistics, but it does
so in a mostly automatic way.
\end{example}

\section{Concluding Remarks}
\label{S:discuss}

Professor Mayo is to be congratulated for her contribution. Besides
resolving the controversy surrounding Birnbaum's theorem, her paper is
an invitation for a fresh discussion on the foundations of statistical
inference. Though LP no longer constrains the frequentist approach, we
have argued here that something more than the basic sampling model is
required for valid statistical inference. The IM framework features the
prediction of unobserved auxiliary variables as this ``something
more,'' and the idea of reducing the dimension of the auxiliary
variable before prediction leads to improved efficiency, accomplishing
what SP and CP are meant to do. We expect further developments for and
from IMs in years to come.

\section*{Acknowledgments}
Supported in part by  NSF Grants DMS-10-07678,
DMS-12-08833 and DMS-12-08841.



\end{document}